\documentclass[12pt,reqno]{amsart}

\usepackage[T1]{fontenc}
\usepackage{amsmath,amssymb,mathtools}
\usepackage{stix2}
\usepackage{geometry}
\usepackage{enumitem}
\usepackage{indentfirst}
\usepackage{microtype}
\usepackage{fancyhdr}
\usepackage[colorlinks,linkcolor=red,citecolor=blue,urlcolor=blue]{hyperref}
\hypersetup{linktocpage=true,
  pdftitle={Optimal first eigenvalue estimates for positive holomorphic sectional curvature}}

\theoremstyle{plain}
\newtheorem{thm}{Theorem}[section]
\newtheorem{lemma}[thm]{Lemma}
\newtheorem{prop}[thm]{Proposition}
\newtheorem{cor}[thm]{Corollary}
\theoremstyle{definition}
\newtheorem{example}[thm]{Example}
\newtheorem{remark}[thm]{Remark}

\newcommand{\btheorem}{\begin{thm}}
\newcommand{\etheorem}{\end{thm}}
\newcommand{\blemma}{\begin{lemma}}
\newcommand{\elemma}{\end{lemma}}
\newcommand{\bproposition}{\begin{prop}}
\newcommand{\eproposition}{\end{prop}}
\newcommand{\bcorollary}{\begin{cor}}
\newcommand{\ecorollary}{\end{cor}}
\newcommand{\bexample}{\begin{example}}
\newcommand{\eexample}{\end{example}}
\newcommand{\bremark}{\begin{remark}}
\newcommand{\eremark}{\end{remark}}
\newcommand{\bproof}{\begin{proof}}
\newcommand{\eproof}{\end{proof}}
\newcommand{\beq}{\begin{equation}}
\newcommand{\eeq}{\end{equation}}
\newcommand{\be}{\begin{eqnarray*}}
\newcommand{\ee}{\end{eqnarray*}}

\newcommand{\p}{\partial}
\newcommand{\bp}{\overline{\partial}}
\renewcommand{\bar}{\overline}
\renewcommand{\phi}{\varphi}
\newcommand{\eps}{\varepsilon}
\newcommand{\C}{\mathbb C}
\newcommand{\R}{\mathbb R}

\renewcommand{\P}{\mathbb P}
\newcommand{\cD}{\mathcal D}
\newcommand{\cR}{\mathcal R}
\newcommand{\FS}{\mathrm{FS}}
\DeclareMathOperator{\Div}{div}
\DeclareMathOperator{\Ric}{Ric}
\DeclareMathOperator{\HSC}{HSC}
\DeclareMathOperator{\tr}{tr}
\DeclareMathOperator{\Vol}{Vol}
\newcommand{\la}{\langle}
\newcommand{\ra}{\rangle}
\newcommand{\dmu}{\,d\mu}
\newcommand{\nm}[1]{\left\|#1\right\|}
\numberwithin{equation}{section}

\setlist[enumerate]{label=(\arabic*),leftmargin=*}

\title{Optimal first eigenvalue estimates for positive holomorphic sectional curvature}

\author{Jiaxuan Fan}
\author{Mingwei Wang}
\author{Xiaokui Yang}
\date{}

\address{Jiaxuan Fan,  Qiuzhen College, Tsinghua University, Beijing, 100084, China}
\email{jxfan26@mails.tsinghua.edu.cn}

\address{Mingwei Wang,  Qiuzhen College, Tsinghua University, Beijing, 100084, China}
\email{wangmw26@mails.tsinghua.edu.cn}

\address{Xiaokui Yang, Department of Mathematics and Yau Mathematical Sciences Center, Tsinghua University, Beijing, 100084, China}
\email{xkyang@mail.tsinghua.edu.cn}
\begin{document}

\begin{abstract}
Let $(M,\omega_g)$ be a compact K\"ahler manifold with
holomorphic sectional curvature  $\geq 2$. We prove that
the first positive eigenvalue of the Laplacian $\Delta_d$ satisfies
$$\lambda_1\geq 4, $$ and $\lambda_1=4$ if and only if $(M,\omega_g)$ is holomorphically
isometric to $(\C\P^1,\omega_{\FS})$.
\end{abstract}

\maketitle

\section{Introduction}

Lower bounds for the first eigenvalue provide a basic connection between
curvature and analysis on a compact manifold. The classical theorem of
Lichnerowicz \cite{Lic58} asserts that a connected closed Riemannian
manifold of real dimension $m\geq 2$ with
$\Ric(g)\geq (m-1)g$ satisfies $\lambda_1\geq m$. Obata
\cite{Oba62} proved that equality characterizes the unit round sphere.
In K\"ahler geometry, the corresponding estimate under
$\Ric(g)\geq(n+1)g$ is $\lambda_1\geq 2(n+1)$. We refer to
\cite{Fut88}, \cite{CWZ25}, \cite{WY26} and the references therein. These results suggest asking for the optimal spectral
bound under a lower bound for holomorphic sectional curvature.\\

Throughout this paper, $(M,\omega_g)$ is a connected compact K\"ahler
manifold without boundary, of complex dimension $n\geq 1$. We write
$\omega=\omega_g$ and $d\mu=\omega^n/n!$. Our Laplacian is the
nonnegative Hodge Laplacian
\beq\label{laplacian-intro}
\Delta_d=dd^*+d^*d,
\qquad \Delta_d=2\Delta_{\bp}\quad\text{on functions}.
\eeq
The first positive eigenvalue of $\Delta_d$ is denoted by $\lambda_1$.
For $0\neq V\in T^{1,0}M$, we define
\beq\label{hsc-intro}
\HSC(V):=\frac{R(V,\bar V,V,\bar V)}{|V|^4}.
\eeq
The Fubini--Study metric $\omega_{\FS}$ on $\P^n$ is normalized by
$\HSC(\omega_{\FS})\equiv 2$. In particular,
$(\C\P^1,\omega_{\FS})$ is the round sphere of Gaussian curvature
$2$, and its first eigenvalue is $4$.\\

The second and third  authors proved in  \cite[Theorem 1.4]{WY26} that $\HSC\geq 2$
implies
\beq\label{previous-bound}
\lambda_1\geq\frac{320(n-1)+576}{81(n-1)+144}.
\eeq
The main result of this paper obtains the optimal first eigenvalue estimates in this setting and confirms the conjecture in
\cite{WY26} affirmatively:

\btheorem\label{MainTheoremA}
Let $(M,\omega_g)$ be a  compact K\"ahler manifold. If $\HSC\geq 2$, then
\beq\label{main-estimate}
\lambda_1\geq 4,
\eeq and
the equality holds if and only if $(M,\omega_g)$ is holomorphically
isometric to $(\C\P^1,\omega_{\FS})$.
\etheorem

\noindent 
For $n\geq 2$ and $ 0<t<\frac{1}{n}$,  consider the $n$-dimensional compact K\"ahler manifold
$$
M_t=
\bigl(\mathbb{CP}^1,(1-(n-1)t)\omega_{\mathrm{FS}}\bigr)
\times\prod_{j=2}^{n}
\bigl(\mathbb{CP}^1,t\omega_{\mathrm{FS}}\bigr).
$$
A straightforward computation shows that  $\mathrm{HSC}\geq 2$
 and $$
\lambda_1(M_t)
=\frac4{1-(n-1)t}\rightarrow4
$$ as $t\rightarrow 0$. Therefore, the lower bound $\lambda_1\geq 4$ is optimal in every dimension.\\

\noindent\textbf{Acknowledgements}.  
 The third author is particularly grateful to Bing-Long Chen and Valentino Tosatti for their insightful comments and suggestions. The third  author is partially supported by National Key R\&D
Program of China 2022YFA1005400 and NSFC grants (No. 12325103, No. 12171262 and No. 12141101).

\vskip 1\baselineskip

\section{Background materials and lemmas}\label{sec-background}
Let $(M,\omega_g)$ be a compact K\"ahler manifold of complex dimension
$n$, and let $\nabla$ denote its Levi-Civita connection. For a real
smooth function $f$, set $\phi=\p f$ and let
$V=(\nabla f)^{1,0}$ be the vector dual to $\bar\phi$. In local
holomorphic coordinates,
\beq\label{complex-gradient}
V^i=g^{i\bar j}f_{\bar j},\qquad
\nabla f=V+\bar V,\qquad
|\nabla f|^2=2|V|^2=2|\phi|^2.
\eeq
Here and below, all $L^p$ norms and integrals are taken with respect to
$d\mu=\omega_g^n/n!$.
For a smooth $(1,0)$ vector field $F$, define
\[
\cD F:=\operatorname{Re}(\nabla_iF^i).
\]
Applying the divergence theorem to $f(F+\bar F)$ therefore gives
\beq\label{integration-parts}
\int_M f\,\cD F\dmu
=-\operatorname{Re}\int_M\phi_iF^i\dmu.
\eeq
In particular, if $u$ is a real smooth function and
$W=(\nabla u)^{1,0}$, then
\[ \cD W=
\nabla_iW^i
=g^{i\bar j}\nabla_i\nabla_{\bar j}u
=-\Delta_{\bp}u
\]
consistently with \eqref{laplacian-intro}.
Let $\eta$ be a smooth closed real one-form and let
$Z=\eta^\sharp$ be its metric dual. Define the smooth endomorphism
$A=A_\eta$ of $T_\R M$ by
\beq\label{real-notation}
A(X)=\nabla_XZ.
\eeq
Metric compatibility and the torsion-free property of $\nabla$ give
\[
(d\eta)(X,Y)
=g(\nabla_XZ,Y)-g(\nabla_YZ,X).
\]
Thus $d\eta=0$ implies that $A$ is symmetric. We also write
$A$ as:
\[
A(X,Y):=g(AX,Y).
\]
Set
\beq\label{hermitian-hessian}
B:=\frac12(A-JAJ).
\eeq
The endomorphism $B$ is symmetric and commutes with $J$ and
$
g(A,B)=|B|^2.
$
Define the real vector fields
\beq\label{currents-definition}
K:=A(JZ,JZ)Z-A(Z,JZ)JZ,
\qquad
K_0:=(\tr A)Z-2BZ.
\eeq
We also write
\[
W=Z^{1,0}=\frac12(Z-\sqrt{-1}JZ),
\qquad
Z=W+\bar W,
\qquad
JZ=\sqrt{-1}(W-\bar W).
\]
With the curvature convention
$
R(X,Y)=[\nabla_X,\nabla_Y]-\nabla_{[X,Y]}$
we have
\beq\label{real-complex-curvature}
\cR:=g(R(Z,JZ)JZ,Z)
=4R(W,\bar W,W,\bar W).
\eeq

\blemma\label{divergence-lemma}
With the notation above,
\begin{align}
	\Div K
	&=-\cR+A(JZ,JZ)\tr A
	+3\la J(AZ),A(JZ)\ra-|A(JZ)|^2,
	\label{div-K}\\
	\Div K_0
	&=(\tr A)^2-2|B|^2.
	\label{div-Kzero}
\end{align}
\elemma

\bproof
Since $A$ is symmetric and $\nabla J=0$, we have
\[
\Div Z=\tr A,
\qquad
\Div(JZ)=\tr(JA)=0.
\]
Consequently,
\[
\Div K
=A(JZ,JZ)\tr A
+Z\bigl(A(JZ,JZ)\bigr)
-(JZ)\bigl(A(Z,JZ)\bigr).
\]
Using $AX=\nabla_XZ$, metric compatibility, and $\nabla J=0$, we obtain
\begin{align*}
	Z\bigl(A(JZ,JZ)\bigr)
	-(JZ)\bigl(A(Z,JZ)\bigr)
	&=\la
	(\nabla_Z\nabla_{JZ}-\nabla_{JZ}\nabla_Z)Z,JZ
	\ra
	+2\la J(AZ),A(JZ)\ra\\
	&=-\cR+\la A([Z,JZ]),JZ\ra
	+2\la J(AZ),A(JZ)\ra.
\end{align*}
Here we used the curvature commutator and
$\la R(Z,JZ)Z,JZ\ra=-\cR$. Moreover,
\[
[Z,JZ]
=\nabla_Z(JZ)-\nabla_{JZ}Z
=J(AZ)-A(JZ).
\]
The symmetry of $A$ therefore gives
\[
\la A([Z,JZ]),JZ\ra
=\la[Z,JZ],A(JZ)\ra
=\la J(AZ),A(JZ)\ra-|A(JZ)|^2.
\]
Combining these identities proves \eqref{div-K}. For the second identity, we use the formula
\[
\Div(\nabla_XX)
=X(\Div X)+\tr\bigl((\nabla X)^2\bigr)+\Ric(X,X),
\]
 for a real vector field $X$.
Since $\nabla J=0$,
\[
2BZ=(A-JAJ)Z
=\nabla_ZZ-\nabla_{JZ}(JZ),
\]
and hence
\[
K_0=(\Div Z)Z-\nabla_ZZ+\nabla_{JZ}(JZ).
\]
Applying the preceding divergence formula, and using
\[
\Div(JZ)=0,
\qquad
\nabla(JZ)=JA,
\qquad
\Ric(JZ,JZ)=\Ric(Z,Z),
\]
we obtain
\begin{align*}
	\Div K_0
	={}&Z(\Div Z)+(\Div Z)^2
	-Z(\Div Z)-\tr(A^2)-\Ric(Z,Z)\\
	&+\tr(JAJA)+\Ric(JZ,JZ)\\
	={}&(\tr A)^2-\tr(A^2)+\tr(JAJA).
\end{align*}
Finally, since $B$ is symmetric and $J^2=-I$,
\[
2|B|^2
=\frac12\tr\bigl((A-JAJ)^2\bigr)
=\tr(A^2)-\tr(JAJA).
\]
Therefore,
$
\Div K_0=(\tr A)^2-2|B|^2
$ and \eqref{div-Kzero} holds.
\eproof

\vskip 1\baselineskip

\section{A weighted divergence inequality}\label{sec-weighted}

The following proposition is the central estimate. 

\bproposition\label{weighted-proposition}
Let $(M,\omega_g)$ be a compact K\"ahler manifold and let $\eta$ be
a smooth closed real one-form. If $W=(\eta^\sharp)^{1,0}$ and
$q=|W|^2$, then
\beq\label{weighted-inequality}
\int_M \left(\cD(qW)\right)^2\dmu
\geq\int_M q\,R(W,\bar W,W,\bar W)\dmu.
\eeq
\eproposition
\bproof
Let $Z=\eta^\sharp$ and retain the notation $A,B,K,K_0$ from
Lemma \ref{divergence-lemma}. Set
\[
s=|Z|^2,\qquad t=\tr A,\qquad
C:=A-B=\frac12(A+JAJ),
\]
and write
\[
b:=\la BZ,Z\ra,\qquad c:=\la CZ,Z\ra.
\] Here $\la,\ra$ is the point-wise inner product.
Since $BJ=JB$ and $CJ=-JC$, we have
\[
A(Z,Z)=b+c,\qquad A(JZ,JZ)=b-c,
\]
as well as
\[
2BZ=AZ-JA(JZ),\qquad
2CZ=AZ+JA(JZ).
\]
 Since $\nabla s=2AZ$,
the product rule gives
\begin{align*}
	\Div(sZ)
	&=st+2(b+c), \qquad
	\la\nabla s,K\ra
	=2\bigl(b^2-c^2-A(Z,JZ)^2\bigr),\\
	\la\nabla(s^2),K_0\ra
	&=4s\bigl(t(b+c)-2\la AZ,BZ\ra\bigr).
\end{align*}
Substituting these expressions into \eqref{div-K} and
\eqref{div-Kzero}, we obtain
$$
	\Div\left(2sK+\frac12s^2K_0\right)+2s\cR=
	\frac12s^2t^2+4stb+4b^2-4c^2
	-s^2|B|^2-8s|CZ|^2-4A(Z,JZ)^2.
$$
Here the terms involving $AZ$ and $A(JZ)$ simplify by
\[
6\la J(AZ),A(JZ)\ra
-2|A(JZ)|^2-4\la AZ,BZ\ra
=-8|CZ|^2.
\]
Consequently,
$$
	[\Div(sZ)]^2-2s\cR
	-\Div\left(2sK+\frac12s^2K_0\right)=
	\frac12(st+4c)^2+s^2|B|^2
	+8s|CZ|^2+8bc+4A(Z,JZ)^2.
$$
The right-hand side is nonnegative by two elementary estimates:
\[
s^2|B|^2\geq 2b^2,
\qquad
s|CZ|^2\geq c^2.
\]
The second is the Cauchy--Schwarz inequality. For the first, note that
$B(Z,Z)=B(JZ,JZ)=b$ and that $Z$ and $JZ$ are orthogonal.
When $s>0$, extend $Z/\sqrt s,JZ/\sqrt s$ to a real orthonormal
frame and retain these two diagonal entries in $|B|^2$.
When $s=0$, both estimates are immediate. Thus
\begin{align*}
	&[\Div(sZ)]^2-2s\cR
	-\Div\left(2sK+\frac12s^2K_0\right)\\
	&\qquad\geq
	\frac12(st+4c)^2+2b^2+8c^2+8bc+4A(Z,JZ)^2\\
	&\qquad=
	\frac12(st+4c)^2+2(b+2c)^2+4A(Z,JZ)^2
	\geq0.
\end{align*}
Integrating over the compact manifold and applying the divergence
theorem gives
\[
\int_M[\Div(sZ)]^2\dmu
\geq 2\int_M s\cR\dmu.
\]
Finally, $
s=2q$, $
\cR=4R(W,\bar W,W,\bar W)$ and 
$\Div(sZ)=4\cD(qW)$ gives \eqref{weighted-inequality}.
\eproof

\bcorollary\label{regularized-corollary}
Let $u\in C^\infty(M,\R)$, $W=(\nabla u)^{1,0}$, and
$q=|W|^2$. For every $\eps>0$,
\beq\label{regularized-inequality}
\int_M[\cD((q+\eps)W)]^2\dmu
\geq\int_M(q+\eps)R(W,\bar W,W,\bar W)\dmu.
\eeq
\ecorollary
\bproof
Let $E$ be a flat complex torus and choose a parallel $(1,0)$
vector field $U$ on $E$ with $|U|^2=\eps$. On the product
K\"ahler manifold $\widehat M=M\times E$, set
\[
\widehat W=W+U.
\]
The real one-form dual to $\widehat W+\bar{\widehat W}$ is
$$
\widehat\eta=du+(U+\bar U)^\flat,$$
which is closed because $U$ is parallel. The product metric gives
$
|\widehat W|^2=q+\eps.
$
Since $q$ depends only on $M$ and $U$ is parallel,
\[
\widehat{\cD}\bigl(|\widehat W|^2\widehat W\bigr)
=\cD\bigl((q+\eps)W\bigr).
\]
Moreover, the flatness of $E$ implies
\[
\widehat R(\widehat W,\bar{\widehat W},
\widehat W,\bar{\widehat W})
=R(W,\bar W,W,\bar W).
\]
Apply Proposition \ref{weighted-proposition} to $\widehat\eta$.
Both integrands depend only on $M$, so integrating over $E$
and cancelling the common positive factor gives
\eqref{regularized-inequality}.
\eproof

\vskip 1\baselineskip

\section{Proof of the main theorem}\label{sec-proofs}

On  $(M,\omega_g)$, we choose a real first eigenfunction $f$ normalized by
\beq\label{eigenfunction-normalization}
\Delta_d f=\lambda_1 f,\qquad
\int_M f\dmu=0,\qquad \int_M f^2\dmu=1.
\eeq
 If we set
$\lambda=\frac{\lambda_1}{2}$, $ \phi=\p f$, and
$V=(\nabla f)^{1,0}$, 
integration by parts gives
\beq\label{gradient-energy}
\Delta_{\bp}f=\lambda f,
\qquad \nm{\phi}_{L^2}^2=\nm{V}_{L^2}^2=\lambda.
\eeq

\blemma\label{prescribed-divergence-lemma}
For each $0<\eps\leq 1/2$, there is a unique smooth real function
$u$ with $\int_Mu\dmu=0$ such that,
\beq\label{prescribed-divergence}
\cD((|W|^2+\eps)W)=f.
\eeq
 where
 $W=(\nabla u)^{1,0}$.
\elemma

\bproof
 It is easy to see that the  equation \eqref{prescribed-divergence} is equivalent to
\beq\label{real-prescribed-equation}
\Div\bigl((|\nabla u|^2+2\eps)\nabla u\bigr)=4f.
\eeq
On the space
\[
\mathcal X=\left\{v\in W^{1,4}(M,\mathbb R) \ \big|\ 
\int_Mv\,d\mu=0\right\},
\]
consider
\[
\mathcal J_\eps(v)
=\int_M\left(\frac14|\nabla v|^4
+\eps|\nabla v|^2+4fv\right)\dmu.
\]
Poincar\'e's inequality gives
\[
\mathcal J_\eps(v)
\geq\frac14\nm{\nabla v}_{L^4}^4
-C\nm{\nabla v}_{L^4}.
\]
Thus $\mathcal J_\eps$ is coercive and weakly lower semicontinuous.
It is also strictly convex on $\mathcal X$. The direct
method therefore gives a unique minimizer $u\in\mathcal X$. Its first variation is
\[
\int_M(|\nabla u|^2+2\eps)
\la\nabla u,\nabla v\ra\dmu
=-4\int_Mfv\dmu.
\]
Since $\int_Mf\dmu=0$, this identity holds for every smooth test
function $v$, not only those of mean zero. Hence $u$ solves
\eqref{real-prescribed-equation} weakly. Conversely, every mean-zero
weak solution is a critical point of the strictly convex functional
$\mathcal J_\eps$, and therefore equals $u$.

It remains to justify smoothness. Since
$\Div(\nabla f)=-\lambda_1f$, we can rewrite the equation as
\[
\Div\left((|\nabla u|^2+2\eps)\nabla u
+\frac4{\lambda_1}\nabla f\right)=0.
\]
For fixed $\eps>0$, this equation has smooth coefficients and the
standard growth and ellipticity conditions associated with $p=4$.
Indeed, the derivative of its flux with respect to the gradient
variable $\xi$ is
\[
(|\xi|^2+2\eps)I+2\xi\otimes\xi.
\]
The equal-exponent case $p=q=4$ and $\mu=1$ of
\cite[Section 2.6, Corollary~5]{DFM23} gives local
$C^{1,\alpha}$ regularity: apply the theorem on coordinate balls
with the boundary trace of $u$, and use strict monotonicity to
identify the resulting solution with $u$.
Thus $u\in C^{1,\alpha}_{\mathrm{loc}}(M)$. Since $\eps>0$,
the linearized coefficients are uniformly elliptic and H\"older
continuous on smaller coordinate balls. Difference quotients
and Schauder estimates give $u\in C^{2,\alpha}_{\mathrm{loc}}(M)$,
and bootstrapping yields $u\in C^\infty(M)$.
This proves the lemma.
\eproof

\vskip 1\baselineskip

\bproof[Proof of  Theorem \ref{MainTheoremA}]
Throughout this proof, $\la\cdot,\cdot\ra_{L^2}$ denotes the
Hermitian $L^2$ inner product induced by $g$.
Recall that $f$ is a real first eigenfunction normalized by
$\nm{f}_{L^2}=1$. With $\lambda=\lambda_1/2$, we have
\[
\Delta_{\bp}f=\lambda f,
\qquad
\nm{\p f}_{L^2}^2=\nm{V}_{L^2}^2=\lambda.
\]
Let $u_\eps$ be the solution from
Lemma \ref{prescribed-divergence-lemma}, and write
\[
W_\eps=(\nabla u_\eps)^{1,0},
\qquad q_\eps=|W_\eps|^2.
\]
Set
\beq\label{two-fluxes}
Y_\eps=q_\eps W_\eps,
\qquad
F_\eps=(q_\eps+\eps)W_\eps.
\eeq
By Corollary \ref{regularized-corollary} and $\HSC\geq2$,
\beq\label{flux-bound}
\begin{aligned}
	2\nm{Y_\eps}_{L^2}^2
	=2\int_M q_\eps^3\dmu
	&\leq\int_M q_\eps R(W_\eps,\bar W_\eps,W_\eps,\bar W_\eps)\dmu\\
	&\leq\int_M(q_\eps+\eps)R(W_\eps,\bar W_\eps,W_\eps,\bar W_\eps)\dmu\\
	&\leq\int_Mf^2\dmu=1.
\end{aligned}
\eeq
H\"older's inequality consequently gives
\beq\label{flux-error}
	\nm{F_\eps-Y_\eps}_{L^2}
	=\eps\nm{W_\eps}_{L^2}
	\leq\eps\Vol(M)^{1/3}
	\left(\int_M q_\eps^3\dmu\right)^{1/6}
	\leq C\eps,
\eeq
where $C=2^{-1/6}\Vol(M)^{1/3}$ is independent of $\eps$.
Thus
\beq\label{flux-limsup}
\nm{F_\eps}_{L^2}\leq\frac1{\sqrt2}+C\eps.
\eeq
Since $\cD F_\eps=f$, integration by parts and the
Cauchy--Schwarz inequality imply
\beq\label{eigenfunction-pairing}
\begin{aligned}
1
&=\la\cD F_\eps,f\ra_{L^2}
=-\operatorname{Re}\la F_\eps,V\ra_{L^2}\\
&\leq\nm{V}_{L^2}\nm{F_\eps}_{L^2}
\leq\sqrt\lambda\left(\frac1{\sqrt2}+C\eps\right).
\end{aligned}
\eeq
Letting $\eps\downarrow0$ yields $\lambda\geq2$, and hence
$
\lambda_1=2\lambda\geq4.
$

\vskip 1\baselineskip

Suppose that $\lambda_1=4$, so that $\lambda=2$. Retain the
normalization of $f$ and the notation $V,Y_\eps,F_\eps$ from the
preceding proof. By \eqref{eigenfunction-pairing} and
\eqref{flux-limsup},
\[
\frac1{\sqrt2}\leq\nm{F_\eps}_{L^2}
\leq\frac1{\sqrt2}+C\eps.
\]
Since $V$ is dual to $\bar\phi$, the exact pairing in
\eqref{eigenfunction-pairing} and $\nm{V}_{L^2}^2=2$ give
\beq\label{equality-strong-limit}
\begin{aligned}
\nm{F_\eps+V/2}_{L^2}^2
&=\nm{F_\eps}_{L^2}^2+\frac14\nm{V}_{L^2}^2
  +\operatorname{Re}\la F_\eps,V\ra_{L^2}\\
&=\nm{F_\eps}_{L^2}^2-\frac12\rightarrow0.
\end{aligned}
\eeq
Together with \eqref{flux-error}, this implies
$Y_\eps\rightarrow-V/2$ strongly in $L^2$.
We next show that the curvature bound is attained in the direction
of $V$. Define the function on $T^{1,0}M$ by
\[
\mathcal H(T)=
\begin{cases}
\displaystyle\frac{R(T,\bar T,T,\bar T)}{|T|^2}-2|T|^2,
&T\neq0,\\[4pt]
0,&T=0.
\end{cases}
\]
This function is continuous, including at the zero section, since
the curvature tensor is bounded on $M$. For $Y_\eps=q_\eps W_\eps$,
\eqref{flux-bound} gives
\beq\label{equality-curvature-excess}
\begin{aligned}
0\leq\int_M\mathcal H(Y_\eps)\dmu
&=\int_M q_\eps R(W_\eps,\bar W_\eps,W_\eps,\bar W_\eps)\dmu
  -2\nm{Y_\eps}_{L^2}^2\\
&\leq1-2\nm{Y_\eps}_{L^2}^2\rightarrow0.
\end{aligned}
\eeq
Choose a sequence $\eps\downarrow0$ along which
$Y_\eps\to-V/2$ almost everywhere. Fatou's lemma yields
$$\int_M\mathcal H(-V/2)\dmu=0.$$ Since the integrand is
nonnegative and continuous, it vanishes everywhere. Thus
\beq\label{equality-directional-curvature}
R(V,\bar V,V,\bar V)=2|V|^4.
\eeq
By \cite[Lemma~6.1]{Yan18}, one deduces that 
\beq \label{equality-directional-ricci}
\Ric(V,\bar V)
\geq(n+1)|V|^2.
\eeq
Let $s=\p f \in \Gamma(M,E)$ with $E=T^{*1,0}M$.  The Bochner-Kodaira formula gives
$$\Delta_{\bp_E } s=\Delta_{\p _E}s+[\sqrt{-1} R^E,\Lambda_\omega]s.$$
By pairing with $s$, we have 
\beq\label{kahler-bochner}
\la\Delta_{\bp_E}s,s\ra_{L^2}
\geq \int_M
\Ric(V,\bar V)\dmu.
\eeq
Since $(M,\omega_g)$ is K\"ahler, one has 
$$\Delta_{\bar\partial_E} s=\Delta_{\bar\partial}(\p f)=\p\Delta_{\bar\partial} f=\lambda \p f.$$
Hence, 
$$4=\lambda \nm{\p f}_{L^2}^2 \geq \int_M
\Ric(V,\bar V)\dmu\geq2(n+1).$$
It follows that $n=1$. In complex dimension one, holomorphic
sectional curvature equals Gaussian curvature, so
$\Ric(g)\geq2g$. Since $\lambda_1(g)=4$, the scaled
Lichnerowicz--Obata theorem \cite{Oba62} implies that $(M,g)$
is isometric to the round sphere of radius $1/\sqrt2$.  Hence $(M,\omega_g)$ is holomorphically isometric to
$(\C\P^1,\omega_{\FS})$.
\eproof

\vskip 1\baselineskip

\noindent\textbf{AI Use Disclosure.} The research direction, main techniques, and overall proof strategy were
determined by the authors. Generative AI was used interactively to help
connect ideas, develop intermediate arguments, and improve the exposition.
The authors independently verified the final mathematical arguments and
take full responsibility for the content of the paper.

\end{document}